\documentclass[oneside,english]{amsart}
\usepackage[T1]{fontenc}
\usepackage[latin9]{inputenc}
\usepackage{url}
\usepackage{amsthm}
\usepackage{amssymb}
\usepackage{esint}

\numberwithin{equation}{section} 
\numberwithin{figure}{section} 
\theoremstyle{plain}
\theoremstyle{plain}
\newtheorem{thm}{Theorem}
  \theoremstyle{remark}
  \newtheorem{rem}[thm]{Remark}

\usepackage{babel}

\begin{document}

\title{Computations on Some Hankel Matrices}

\author{Ruiming Zhang}
\begin{abstract}
In this note, we present the determinant, the inverse and a lower
bound for the smallest eigenvalue for some Hankel matrices
\end{abstract}


\curraddr{School of Mathematical Sciences\\
Guangxi Normal University\\
Guilin City, Guangxi 541004\\
P. R. China.}


\email{ruimingzhang@yahoo.com}

\maketitle

\section{Introduction}

For each nonnegative integer $n$, the $n$-th Hilbert matrix is \cite{Weisstein} 

\[
\left(\frac{1}{j+k+1}\right)_{j,k=0}^{n}.\]
These matrices are the moment matrices associated the Legendre polynomials.
The generalized Hilbert matrices, which are also called Hankel matrices,
are from the generalized moment matrices associated with some more
general orthogonal polynomials. Some interesting questions for Hankel
matrices are the determinants, inverses and lower bounds for the smallest
eigenvalues. In \cite{Zhang} we have developed a general method to
compute the determinants, inverses and lower bounds for the smallest
eigenvalues for the generalized moment matrices associated with some
orthogonal systems (not just limited to orthogonal polynomials). In
this note we apply the results to some Hankel matrices. The following
theorem is adapted from \cite{Zhang} and we won't repeat the proof
here.
\begin{thm}
\label{thm:1}Given a probability measure $P(dx)$ on $\mathbb{R}$,
for each nonnegative integer $n$, let\[
\mu_{n}=\int_{\mathbb{R}}x^{n}P(dx),\]

\[
G_{n}=\left(\mu_{j+k}\right)_{j,k=0}^{n}\]
 and\begin{align*}
p_{n}(x) & =\sum_{k=0}^{n}a_{n,k}x^{k},\quad n=0,1,\dots\end{align*}
be the orthonormal polynomials, then 

\begin{align*}
\det G_{n} & =\prod_{j=0}^{n}a_{j,j}^{-2}\end{align*}
and\begin{align*}
G_{n}^{-1} & =\left(\gamma_{j,k}\right)_{j,k=0}^{n}\end{align*}
with\begin{align*}
\gamma_{j,k} & =\sum_{\ell=\max(j,k)}^{n}\overline{a_{\ell,j}}a_{\ell,k}.\end{align*}
Furthermore, if there is a complex number $z_{0}$ with $|z_{0}|=1$
such that for each nonnegative integer $n$, the following sequence
\[
a_{n,k}z_{0}^{k},\quad k=0,1,\dots\]
have the same sign , then the smallest eigenvalue $\lambda_{s}$ of
the matrix $G_{n}$ has a lower bond\begin{align*}
\lambda_{s} & \ge\frac{1}{\sum_{m=0}^{n}|p_{m}(z_{0})|^{2}}.\end{align*}
\end{thm}
\begin{rem}
In the case that all the $p_{m}(z_{0})$ are real, we could apply
the Christoffel-Darboux formula to get \cite{Andrews,Szego} \[
\lambda_{s}\ge\frac{a_{n+1,n+1}}{a_{n,n}\left\{ p'_{n+1}(z_{0})p_{n}(z_{0})-p_{n+1}(z_{0})p'_{n}(z_{0})\right\} }.\]
 
\end{rem}
Recall that the Euler's $\Gamma(z)$ is defined as \cite{Andrews,Koekoek,Szego}
\begin{align*}
\Gamma(z) & =\int_{0}^{\infty}x^{z-1}e^{-x}dx,\quad\Re(z)>0\end{align*}
and it could be analytically extended to a meromorphic function on
the complex plane. The shifted factorial of $z$ is defined as\begin{align*}
(z)_{n} & =\frac{\Gamma(z+n)}{\Gamma(z)},\quad n\in\mathbb{Z}.\end{align*}
The hypergeometric function ${}_{2}F_{1}$ is defined as\begin{align*}
{}_{2}F_{1}\left(\begin{array}{c}
a,b\\
c\end{array};z\right) & =\sum_{n=0}^{\infty}\frac{(a)_{n}(b)_{n}}{(c)_{n}n!}z^{n}\end{align*}
for $|z|<1$. Euler's Beta integral could be evaluated in terms of
$\Gamma(z)$,\begin{align*}
\int_{0}^{1}x^{\alpha-1}(1-x)^{\beta-1}dx & =\frac{\Gamma(\alpha)\Gamma(\beta)}{\Gamma(\alpha+\beta)},\quad\Re(\alpha),\Re(\beta)>0.\end{align*}
For any complex number $a$ and $0<q<1$, we define \cite{Andrews,Koekoek}

\begin{align*}
(a;q)_{\infty} & =\prod_{m=0}^{\infty}(1-aq^{m}),\quad(a;q)_{m}=\frac{(a;q)_{\infty}}{(aq^{m};q)_{\infty}}.\end{align*}
The $q$-Binomial theorem is\begin{align*}
\frac{(az;q)_{\infty}}{(z;q)_{\infty}} & =\sum_{k=0}^{\infty}\frac{(a;q)_{k}}{(q;q)_{k}}z^{k},\quad|z|<1,\end{align*}
one of its direct consequences is \begin{align*}
(z;q)_{\infty} & =\sum_{k=0}^{\infty}\frac{q^{\binom{k}{2}}\left(-z\right)^{k}}{(q;q)_{k}}.\end{align*}
The confluent $q$-hypergeometric series ${}_{1}\phi_{1}$ \begin{align*}
{}_{1}\phi_{1}\left(a;b;q,z\right) & =\sum_{n=0}^{\infty}\frac{(a;q)_{n}q^{\binom{n}{2}}(-z)^{n}}{(b;q)_{n}(q;q)_{n}}\end{align*}
satisfies the following identity,\begin{align*}
{}_{1}\phi_{1}\left(a;b;q,b/a\right) & =\frac{(b/a;q)_{\infty}}{(b;q)_{\infty}}.\end{align*}
 Hold $|b|<1$ fixed and let $a\to\infty$ in the above formula to
obtain the Cauchy's formula \[
\sum_{n=0}^{\infty}\frac{q^{n(n-1)}b^{n}}{(q,b;q)_{n}}=\frac{1}{(b;q)_{\infty}}.\]

\section{Applications}

\subsection{Laguerre Polynomials}

The Laguerre polynomials $\left\{ L_{n}^{\alpha}(x)\right\} _{n=0}^{\infty}$
are defined as \cite{Andrews,Koekoek,Szego} \begin{align*}
L_{n}^{\alpha}(x) & =\frac{(\alpha+1)_{n}}{n!}\sum_{k=0}^{n}\frac{(-n)_{k}x^{k}}{(\alpha+1)_{k}k!}\end{align*}
 for $n\ge0$ and we assume that \begin{align*}
L_{-1}^{\alpha}(x) & =0.\end{align*}
For any $\alpha>-1$, the orthogonal relation for the Laguerre polynomials
is

\begin{align*}
\int_{0}^{\infty}L_{m}^{\alpha}(x)L_{n}^{\alpha}(x)\frac{x^{\alpha}e^{-x}}{\Gamma(\alpha+1)}dx & =\frac{(\alpha+1)_{n}}{n!}\delta_{mn}\end{align*}
for any nonnegative integers $m,n$ where\begin{align*}
\delta_{mn} & =\begin{cases}
1 & m=n\\
0 & m\neq n\end{cases}.\end{align*}
Clearly, the $n$-th moment is\begin{align*}
m_{n} & =\frac{\int_{0}^{\infty}x^{\alpha+n}e^{-x}dx}{\Gamma(\alpha+1)}=(\alpha+1)_{n},\end{align*}
and

\begin{align*}
\ell_{n}^{(\alpha)}(x) & =(-1)^{n}\sqrt{\frac{n!}{(\alpha+1)_{n}}}L_{n}^{\alpha}(x)\end{align*}
is the $n$-th orthonormal polynomial. Apply Theorem \ref{thm:1}
with\begin{align*}
a_{n,k} & =\sqrt{\frac{(\alpha+1)_{n}}{n!}}\frac{(-n)_{k}(-1)^{n}}{(\alpha+1)_{k}k!}\end{align*}
we get \begin{align*}
\det\left((\alpha+1)_{j+k}\right)_{0\le j,k\le n} & =\prod_{k=0}^{n}\left\{ k!(\alpha+1)_{k}\right\} ,\end{align*}
and its inverse matrix is\begin{align*}
\left((\alpha+1)_{j+k}\right)_{0\le j,k\le n}^{-1} & =\left(\sum_{\ell=0}^{n}\frac{(\alpha+1)_{\ell}(-\ell)_{j}(-\ell)_{k}}{\ell!(\alpha+1)_{j}(\alpha+1)_{k}j!k!}\right)_{j,k=0}^{n}.\end{align*}
For $\alpha>-1$, the smallest eigenvalue of the matrix $\left((\alpha+1)_{j+k}\right)_{0\le j,k\le n}$
is \begin{align*}
\lambda_{s} & \ge\left\{ \sum_{\ell=0}^{n}\frac{\ell!}{(\alpha+1)_{\ell}}L_{\ell}^{(\alpha)}(-1)^{2}\right\} ^{-1}\\
 & =\frac{(\alpha+1)_{n}}{(n+1)!}\frac{1}{L_{n}^{(\alpha+1)}(-1)L_{n}^{(\alpha)}(-1)-L_{n+1}^{(\alpha)}(-1)L_{n-1}^{(\alpha+1)}(-1)},\end{align*}
where we applied the formula\[
\frac{d}{dx}L_{n}^{(\alpha)}(x)=-L_{n-1}^{(\alpha+1)}(x).\]
From the the Perron formula \cite{Szego} \begin{align*}
L_{n}^{(\alpha)}(x) & =\frac{e^{x/2}}{2\pi}(-x)^{-\alpha/2-1/4}n^{\alpha/2-1/4}\exp\left\{ 2(-nx)^{1/2}\right\} \left\{ 1+\mathcal{O}\left(\frac{1}{n^{1/2}}\right)\right\} \end{align*}
for $x\in\mathbb{C}\backslash(0,\infty)$ as $n\to\infty$ to get\begin{align*}
L_{n}^{(\alpha+1)}(-1)L_{n}^{(\alpha)}(-1)-L_{n+1}^{(\alpha)}(-1)L_{n-1}^{(\alpha+1)}(-1) & =\frac{n^{\alpha-1}\exp(4n^{1/2})}{8\pi e}\left\{ 1+\mathcal{O}\left(\frac{1}{n^{1/2}}\right)\right\} \end{align*}
as $n\to\infty$, this together with \begin{align*}
\frac{\Gamma(n+\alpha)}{\Gamma(n+\beta)} & =n^{\alpha-\beta}\left\{ 1+\mathcal{O}\left(\frac{1}{n}\right)\right\} \end{align*}
as $n\to\infty$ gives\begin{align*}
\left\{ \sum_{\ell=0}^{n}\frac{\ell!}{(\alpha+1)_{\ell}}L_{\ell}^{(\alpha)}(-1)^{2}\right\} ^{-1} & =\frac{8\pi e}{\Gamma(\alpha+1)}\exp(-4n^{1/2})\left\{ 1+\mathcal{O}\left(\frac{1}{n^{1/2}}\right)\right\} \end{align*}
as $n\to\infty$.

\subsection{The Jacobi Polynomials $\{P_{n}^{(\alpha,\beta)}(x)\}_{n=0}^{\infty}$}

The Jacobi polynomials $\left\{ P_{n}^{(\alpha,\beta)}(x)\right\} _{n=0}^{\infty}$
are defined as \cite{Andrews,Koekoek,Szego}\begin{align*}
P_{n}^{(\alpha,\beta)}(x) & =\frac{(\alpha+1)_{n}}{n!}{}_{2}F_{1}\left(\begin{array}{c}
-n,n+\alpha+\beta+1\\
\alpha+1\end{array};\frac{1-x}{2}\right)\end{align*}
for $n\ge0$ and \begin{align*}
P_{-1}^{(\alpha,\beta)}(x) & =0.\end{align*}
For $\alpha,\beta>-1$, they satisfy the orthogonal relation

\begin{align*}
\int_{-1}^{1}P_{m}^{(\alpha,\beta)}(x)P_{n}^{(\alpha,\beta)}(x)w(x)dx & =h_{n}\delta_{mn}\end{align*}
for all nonnegative integers $n,m$ where\begin{align*}
w(x) & =(1-x)^{\alpha}(1+x)^{\beta},\end{align*}
 and\begin{align*}
h_{n} & =\frac{2^{\alpha+\beta+1}\Gamma(\alpha+n+1)\Gamma(\beta+n+1)}{(2n+\alpha+\beta+1)\Gamma(\alpha+\beta+n+1)n!}.\end{align*}
Since \begin{align*}
P_{n}^{(\alpha,\beta)}(x) & =\frac{(\beta+1)_{n}}{(-1)^{n}n!}{}_{2}F_{1}\left(\begin{array}{c}
-n;n+\alpha+\beta+1\\
\beta+1\end{array};\frac{1+x}{2}\right),\end{align*}
we let\begin{align*}
A_{n}^{(\alpha,\beta)}(y) & =P_{n}^{(\alpha,\beta)}(2y-1)\end{align*}
then, for $\alpha,\beta>-1$ we have

\begin{align*}
\int_{0}^{1}A_{m}^{(\alpha,\beta)}(y)A_{n}^{(\alpha,\beta)}(y)\tilde{w}(y)dy & =h_{n}^{(\alpha,\beta)}\delta_{mn}\end{align*}
for all nonnegative integers $n,m$ where\begin{align*}
\tilde{w}(y) & =\frac{y^{\alpha}(1-y)^{\beta}\Gamma(\alpha+\beta+2)}{\Gamma(\alpha+1)\Gamma(\beta+1)},\end{align*}
and\[
h_{n}^{(\alpha,\beta)}=\frac{(\alpha+\beta+1)(\alpha+1)_{n}(\beta+1)_{n}}{(2n+\alpha+\beta+1)n!(\alpha+\beta+1)_{n}}.\]
Clearly, the $n$-th moment is \begin{align*}
\mu_{n} & =\int_{0}^{1}y^{n}\tilde{w}(y)dy=\frac{(\alpha+1)_{n}}{(\alpha+\beta+2)_{n}}\end{align*}
and the $n$-th orthonormal polynomial is

\begin{align*}
a_{n}^{(\alpha,\beta)}(y) & =\frac{1}{\sqrt{h_{n}^{(\alpha,\beta)}}}A_{n}^{(\alpha,\beta)}(y),\end{align*}
 or\begin{align*}
a_{n}^{(\alpha,\beta)}(y) & =\sqrt{\frac{(2n+\alpha+\beta+1)(\beta+1)_{n}(\alpha+\beta+1)_{n}}{(\alpha+\beta+1)(\alpha+1)_{n}n!}}\\
\times & (-1)^{n}{}_{2}F_{1}\left(\begin{array}{c}
-n;n+\alpha+\beta+1\\
\beta+1\end{array};y\right).\end{align*}
Thus, \begin{align*}
a_{n,k} & =\sqrt{\frac{(2n+\alpha+\beta+1)(\beta+1)_{n}(\alpha+\beta+1)_{n}}{(\alpha+\beta+1)(\alpha+1)_{n}n!}}\\
\times & \frac{(-n)_{k}(n+\alpha+\beta+1)_{k}}{(-1)^{n}(\beta+1)_{k}k!},\end{align*}
and for each nonnegative integer $n$, we have\begin{align*}
\det\left(\frac{(\alpha+1)_{j+k}}{(\alpha+\beta+2)_{j+k}}\right)_{0\le j,k\le n} & =\prod_{m=0}^{n}\frac{(\beta+1)_{m}(\alpha+1)_{m}m!}{(\alpha+\beta+2)_{2m}(m+\alpha+\beta+1)_{m}}\end{align*}
and\begin{align*}
\left(\frac{(\alpha+1)_{j+k}}{(\alpha+\beta+2)_{j+k}}\right)_{0\le j,k\le n}^{-1} & =\left(\gamma_{j,k}\right)_{0\le j,k\le n}\end{align*}
with\begin{align*}
\gamma_{j,k} & =\sum_{m=0}^{n}\frac{(2m+\alpha+\beta+1)(\beta+1)_{m}(\alpha+\beta+1)_{m}}{(\alpha+\beta+1)(\alpha+1)_{m}m!}\\
 & \frac{(-m)_{j}(-m)_{k}(m+\alpha+\beta+1)_{j}(m+\alpha+\beta+1)_{k}}{j!k!(\beta+1)_{j}(\beta+1)_{k}}.\end{align*}
 Then, its smallest eigenvalue of the matrix $\left(\frac{(\alpha+1)_{j+k}}{(\alpha+\beta+2)_{j+k}}\right)_{0\le j,k\le n}$
has a lower bound

\begin{align*}
\lambda_{s} & \ge\left\{ \sum_{m=0}^{n}\frac{\left(P_{m}^{(\alpha,\beta)}(-3)\right)^{2}}{h_{m}^{(\alpha,\beta)}}\right\} ^{-1}.\end{align*}
 From the Christoffel-Darboux formula and the formula\[
\frac{d}{dx}P_{n}^{(\alpha,\beta)}(x)=\frac{(n+\alpha+\beta+1)}{2}P_{n-1}^{(\alpha+1,\beta+1)}(x),\]
and\[
P_{n}^{(\alpha,\beta)}(-x)=(-1)^{n}P_{n}^{(\beta,\alpha)}(x)\]
to get\begin{align*}
\sum_{m=0}^{n}\frac{\left(P_{m}^{(\alpha,\beta)}(-3)\right)^{2}}{h_{m}^{(\alpha,\beta)}} & =\frac{(\beta+1+n)(\alpha+\beta+1+n)}{(\alpha+\beta+1+2n)(\alpha+\beta+2+2n)h_{n}^{(\alpha,\beta)}}\\
\times & \left\{ (\alpha+\beta+2+n)P_{n}^{(\beta+1,\alpha+1)}(3)P_{n}^{(\beta,\alpha)}(3)\right.\\
- & \left.(\alpha+\beta+1+n)P_{n+1}^{(\beta,\alpha)}(3)P_{n-1}^{(\beta+1,\alpha+1)}(3)\right\} .\end{align*}
 From the asymptotic formula of Jacobi polynomials to obtain \cite{Szego}\[
P_{n}^{(\alpha,\beta)}(3)=\frac{\left(3+2\sqrt{2}\right)^{n}}{n^{1/2}}\left\{ \phi_{0}(\alpha,\beta;3+2\sqrt{2})+\mathcal{O}\left(\frac{1}{n}\right)\right\} \]
 as $n\to\infty$ and\begin{align*}
\sum_{m=0}^{n}\frac{\left(P_{m}^{(\alpha,\beta)}(-3)\right)^{2}}{h_{m}^{(\alpha,\beta)}} & =\phi_{0}(\alpha,\beta;3+2\sqrt{2})\phi_{0}(\alpha+1,\beta+1;3+2\sqrt{2})\\
\times & \frac{\Gamma(\alpha+1)\Gamma(\beta+1)(3+2\sqrt{2})^{2n}}{2\Gamma(\alpha+\beta+2)}\left\{ 1+\mathcal{O}\left(\frac{1}{n}\right)\right\} \end{align*}
 as $n\to\infty$. Thus\[
\left\{ \sum_{m=0}^{n}\frac{\left(P_{m}^{(\alpha,\beta)}(-3)\right)^{2}}{h_{m}^{(\alpha,\beta)}}\right\} ^{-1}=\mathcal{O}\left(\frac{1}{(3+2\sqrt{2})^{2n}}\right)\]
 as $n\to\infty$.

\subsection{The $q$-Laguerre Polynomials $\{L_{n}^{(\alpha)}(x;q)\}_{n=0}^{\infty}$
with $q\in(0,1)$ and $\alpha>-1$}

The $q$-Laguerre polynomials $\{L_{n}^{(\alpha)}(x;q)\}_{n=0}^{\infty}$
are defined as \cite{Andrews,Koekoek} \begin{align*}
L_{n}^{(\alpha)}(x;q) & =\frac{(q^{\alpha+1};q)_{n}}{(q;q)_{n}}\sum_{k=0}^{n}\frac{(q^{-n};q)_{k}q^{\binom{k+1}{2}}q^{(\alpha+n)k}x^{k}}{(q;q)_{k}(q^{\alpha+1};q)_{k}}\end{align*}
for $n\ge0$, and we assume that\[
L_{-1}^{(\alpha)}(x;q)=0.\]
The moment problem of the $q$-Laguerre polynomials is indeterminate
and one of the orthogonality for $\{L_{n}^{(\alpha)}(x;q)\}_{n=0}^{\infty}$
is\begin{align*}
\sum_{k=-\infty}^{\infty}L_{m}^{(\alpha)}(q^{k};q)L_{n}^{(\alpha)}(q^{k};q)w_{ql}(q^{k};\alpha) & =\frac{(q^{\alpha+1};q)_{n}}{(q;q)_{n}q^{n}}\delta_{m,n},\end{align*}
for $m,n\ge0$ with\begin{align*}
w_{ql}(q^{k};\alpha) & =\frac{(q^{\alpha+1};q)_{\infty}(-q;q)_{\infty}(-1;q)_{k}q^{k(\alpha+1)}}{(q;q)_{\infty}(-q^{\alpha+1};q)_{\infty}(-q^{-\alpha};q)_{\infty}}.\end{align*}
Clearly, the $n$-th moment is\begin{align*}
\mu_{n}(\alpha) & =\sum_{k=0}^{\infty}q^{kn}w_{ql}(q^{k};\alpha)=\left(q^{\alpha+1};q\right)_{n}q^{-\alpha n-\binom{n+1}{2}},\end{align*}
and the orthonormal system is given by \begin{align*}
\ell_{n}^{(\alpha)}(x;q) & =\sqrt{\frac{(q^{\alpha+1};q)_{n}q^{n}}{(q;q)_{n}}}\sum_{k=0}^{n}\frac{(q^{-n};q)_{k}q^{\binom{k+1}{2}}q^{(\alpha+n)k}x^{k}}{(-1)^{n}(q;q)_{k}(q^{\alpha+1};q)_{k}},\end{align*}
then,\begin{align*}
a_{n,k} & =\sqrt{\frac{(q^{\alpha+1};q)_{n}q^{n}}{(q;q)_{n}}}\frac{(q^{-n};q)_{k}q^{\binom{k+1}{2}}q^{(\alpha+n)k}}{(-1)^{n}(q;q)_{k}(q^{\alpha+1};q)_{k}}.\end{align*}
According to Theorem \ref{thm:1}, the matrix \begin{equation}
\left((q^{\alpha+1};q)_{j+k}q^{-\binom{j+k+1}{2}-\alpha(j+k)}\right)_{j,k=0}^{n}\label{eq:2.2.1}\end{equation}
 has inverse \[
\left(\gamma_{j,k}\right)_{j,k=0}^{n}\]
 where\begin{align*}
\gamma_{j,k} & =\sum_{m=0}^{n}\frac{(q^{\alpha+1};q)_{m}q^{m(j+k+1)}(q^{-m};q)_{j}(q^{-m};q)_{k}q^{\alpha(j+k)+\binom{j+1}{2}+\binom{k+1}{2}}}{(q;q)_{m}(q;q)_{j}(q;q)_{k}(q^{\alpha+1};q)_{j}(q^{\alpha+1};q)_{k}}.\end{align*}
Its determinant is \begin{align*}
\det\left((q^{\alpha+1};q)_{j+k}q^{-\binom{j+k+1}{2}-\alpha(j+k)}\right)_{j,k=0}^{n} & =\frac{{\displaystyle \prod_{m=0}^{n}(q;q)_{m}(q^{\alpha+1};q)_{m}}}{q^{n(n+1)(4n+6\alpha+5)/6}},\end{align*}
and the smallest eigenvalue has a lower bound is\[
\left(\sum_{m=0}^{n}\left|\ell_{m}^{(\alpha)}(-1;q)\right|^{2}\right)^{-1}.\]
But\begin{align*}
\ell_{m}^{(\alpha)}(-1;q) & =\sqrt{\frac{(q^{\alpha+1};q)_{m}q^{m}}{(q;q)_{m}}}\sum_{k=0}^{m}\frac{(q^{-m};q)_{k}q^{\binom{k}{2}}\left(-q^{m+\alpha+1}\right)^{k}}{(-1)^{m}(q;q)_{k}(q^{\alpha+1};q)_{k}},\\
= & (-1)^{m}\sqrt{\frac{(q^{\alpha+1};q)_{m}q^{m}}{(q;q)_{m}}}{}_{1}\phi_{1}\left(q^{-m};q^{\alpha+1};q,q^{m+\alpha+1}\right)\\
= & (-1)^{m}\sqrt{\frac{(q^{\alpha+1};q)_{m}q^{m}}{(q;q)_{m}}}\frac{(q^{m+\alpha+1};q)_{\infty}}{(q^{\alpha+1};q)_{\infty}}\\
= & (-1)^{m}\sqrt{\frac{q^{m}}{(q,q^{\alpha+1};q)_{m}}},\end{align*}
thus the lower bound for the smallest eigenvalue is\[
\left(\sum_{m=0}^{n}\frac{q^{m}}{(q,q^{\alpha+1};q)_{m}}\right)^{-1}.\]
Notice that everything involved here is a rational function of variable
$q$, we may replace $q$ by $q^{-1}$ then apply the relation\[
(a;q)_{n}=(-a)^{n}q^{\binom{n}{2}}(a^{-1};q^{-1})_{n}\]
and the diagonal matrix \[
D_{n}=\left((-1)^{j}\delta_{j,k}\right)_{j,k=0}^{n},\quad D_{n}^{2}=I\]
to obtain\begin{align*}
\det\left((q^{\alpha+1};q)_{j+k}\right)_{0\le j,k\le n} & =a^{n(n+1)\alpha/2}q^{n(n+1)(2n+1)/6}\prod_{k=0}^{n}(q;q)_{k}(q^{\alpha+1};q)_{k},\end{align*}
 and\begin{align*}
\left((q^{\alpha+1};q)_{j+k}\right)_{0\le j,k\le n}^{-1} & =\left(\sum_{m=0}^{n}\frac{q^{j+k}(q^{\alpha+1};q)_{m}(q^{-m};q)_{j}(q^{-m};q)_{k}}{(q;q)_{j}(q;q)_{k}(q^{\alpha+1};q)_{j}(q^{\alpha+1};q)_{k}(q;q)_{m}q^{(\alpha+1)m}}\right)_{j,k=0}^{n}.\end{align*}
and the smallest eigenvalue has a lower bound\[
\left(\sum_{m=0}^{n}\frac{q^{2\binom{m}{2}}q^{m(\alpha+1)}}{(q,q^{\alpha+1};q)_{m}}\right)^{-1},\]
which is bounded again by the absolute constant\[
\left(\sum_{m=0}^{\infty}\frac{q^{2\binom{m}{2}}q^{m(\alpha+1)}}{(q,q^{\alpha+1};q)_{m}}\right)^{-1}=(q^{\alpha+1};q)_{\infty}.\]

\end{document}